\documentclass[11pt]{article}
\usepackage[left=1in,top=1in,right=1in,nohead,bottom=1in]{geometry}

\usepackage{amsmath}
\usepackage{amssymb}
\usepackage{amsthm}
\usepackage{graphicx}
\usepackage{url}
\usepackage{colordvi}
\usepackage{algorithm}
\usepackage{algorithmic}
\usepackage{verbatim}
\usepackage{color}

\newcommand{\R}{{\mathbb R}}
\newcommand{\Z}{\mathbb{Z}}

\newcommand{\C}{\mathbb{C}}
\newcommand{\F}{\mathbb{F}}
\newcommand{\gen}[1]{\langle #1 \rangle_{\F_2}}

\newcommand{\Q}{\mathbb{Q}}

\newcommand{\set}[1]{\left\{#1\right\}}

\renewcommand{\>}{\rangle}

\renewcommand{\k}{\mathbb K}
\newcommand{\kc}{{\overline{\k}}}

\newtheorem{theorem}{Theorem}[section]
\newtheorem{lemma}[theorem]{Lemma}

\newtheorem{problem}[theorem]{Problem}
\newtheorem{question}[theorem]{Question}
\newtheorem{proposition}[theorem]{Proposition}
\newtheorem{definition}[theorem]{Definition}
\newtheorem{corollary}[theorem]{Corollary}
\newtheorem{remark}[theorem]{Remark}
\newtheorem{example}[theorem]{Example}

\begin{document}
\title{Recognizing Graph Theoretic \\ Properties with Polynomial Ideals}

\author{J.A. De Loera, C. Hillar\footnote{The second author is partially supported by an NSA Young Investigator Grant and an NSF All-Institutes Postdoctoral Fellowship administered by the Mathematical Sciences Research Institute through its core grant DMS-0441170.}, P.N. Malkin, M. Omar \footnote{The fourth author is partially supported by NSERC Postgraduate Scholarship 281174.} \footnote{All other authors are partially
supported by NSF grant DMS-0914107 and an IBM OCR award}}





\maketitle

\abstract{Many hard combinatorial problems can be modeled by a system
  of polynomial equations.  N. Alon coined the term {\em polynomial
    method} to describe the use of nonlinear polynomials when solving
  combinatorial problems.  We continue the exploration of the
  polynomial method and show how the algorithmic theory of polynomial
  ideals can be used to detect $k$-colorability, unique Hamiltonicity,
  and automorphism rigidity of graphs.  Our techniques are diverse and
  involve Nullstellensatz certificates, linear algebra over finite
  fields, Gr\"obner bases, toric algebra, convex programming,
  and real algebraic geometry.}

\section{Introduction}
In his well-known survey \cite{alonsurvey}, Noga Alon used the term {\em polynomial
  method} to refer to the use of nonlinear polynomials when solving combinatorial problems. Although the polynomial method is not yet as  widely used as its linear counterpart, increasing numbers of researchers are using the algebra of multivariate polynomials to solve 
  interesting problems (see for example  \cite{AlonTarsi92,DeLoera95,DeLoeraLeeMalkinMargulies08, Eliahou92,Fischer88, HillarLim, HillarWindfeldt08,Lovasz1994, LiLi81,matiyasevich1,matiyasevich2, Onn,Simis94} and references therein).  In the concluding remarks of \cite{alonsurvey}, Alon asked whether it is possible to modify algebraic proofs to yield efficient algorithmic solutions to combinatorial problems.  In this paper, we explore this question further.  We use polynomial ideals and zero-dimensional varieties to study three hard recognition problems in graph theory.  We show that this approach can be fruitful both theoretically and computationally, and in some cases, result in efficient recognition strategies.

Roughly speaking, our approach is to associate to a combinatorial question (e.g., is a graph $3$-colorable?) a system of polynomial equations $J$ such that the combinatorial problem has a positive answer if and only if system $J$ has a solution.  These highly structured systems of equations (see Propositions \ref{prop:encodings}, \ref{openinghamiltonprop}, and \ref{AutGEncodingprop}), which we refer to as \emph{combinatorial systems of equations}, are then solved using various methods including linear algebra over finite fields, Gr\"obner bases, or semidefinite programming.  As we shall see below this methodology is applicable in a wide range of contexts.

In what follows, $G = (V,E)$ denotes an undirected simple graph on vertex set $V=\set{1,\ldots, n}$ and edges
$E$. Similarly, by $G=(V,A)$ we mean that $G$ is a \emph{directed} graph with arcs $A$. 
When $G$ is undirected, we let \[Arcs(G)=\{(i,j): i,j \in V, \ \text{and} \ \{i,j\} \in E\}\] consist of all possible arcs for each edge in $G$. We study three classical graph problems.

First, in Section \ref{3coloringsec}, we explore $k$-colorability using techniques from commutative algebra and algebraic geometry. 
The following polynomial formulation of $k$-colorability is well-known \cite{Bayer82}. 
\begin{proposition} 
\label{prop:encodings}
Let $G=(V,E)$ be an undirected simple graph on vertices $V = \{1,\ldots,n\}$.
Fix a positive integer $k$, and let $\k$ be a field with characteristic relatively prime to $k$.  The
polynomial system \[J_G = \{x_i^k-1=0, \  x_i^{k-1} + x_i^{k-2}x_j + \cdots + x_j^{k-1} =0 : \ i\in V,  \  \{i,j\}\in E\}\] 
has a common zero over $\overline{\k}$ (the algebraic closure of $\k$) if and only if the graph $G$ is $k$-colorable. 
\end{proposition}

\begin{remark}
Depending on the context, the fields $\k$ we use in this paper will be the rationals $\Q$, the reals $\R$, the complex numbers $\C$, or finite fields $\F_p$ with $p$ a prime number.
\end{remark}

Hilbert's Nullstellensatz \cite[Theorem 2, Chapter 4]{CoxLittleOShea92} states that a
system of polynomial equations $\{f_1(x)=0,\dots,f_r(x)=0\}$ with coefficients in $\k$
has no solution with entries in its algebraic closure $\kc$ if and only if 
\[1=\sum_{i=1}^r \beta_if_i, \ \  \text{ for some polynomials $\beta_1,\ldots,\beta_r \in {\k}[x_1,\dots,x_n]$}.\]  
Thus, if the system has no solution, there is a 
{\em Nullstellensatz certificate} that the associated combinatorial problem is infeasible.
We can find a Nullstellensatz certificate  $1=\sum_{i=1}^r \beta_if_i$ of a
given degree $D:=\max_{1 \leq i \leq r} \{\deg(\beta_i) \}$ or determine
that no such certificate exists by solving a system of \emph{linear equations} whose
variables are in bijection with the coefficients of the monomials of
$\beta_1,\ldots,\beta_r$ (see \cite{ima} and the many references therein).
The number of variables in this linear system grows with the number ${n+D \choose D}$ of
monomials of degree at most $D$.  Crucially, the linear system, which can be thought of as a $D$-th
order linear relaxation of the polynomial system, can be solved in time that is polynomial
in the input size for fixed degree $D$ (see \cite[Theorem 4.1.3]{SusanPhd}  or the survey \cite{ima}).  
The degree $D$ of a Nullstellensatz certificate of an infeasible polynomial system 
cannot be more than known bounds \cite{Kollar1988}, and thus,
by searching for certificates of increasing degrees, we obtain a finite (but potentially long) 
procedure to decide whether a system is feasible or not (this is the NulLA algorithm in 
\cite{SusanPhd,DeLoeraLeeMarguliesOnn08,DeLoeraLeeMalkinMargulies08}).  
The philosophy of ``linearizing'' a system of arbitrary polynomials has also
been applied in other contexts besides combinatorics, including computer algebra \cite{Faugere1999, KehreinKreuzer05,MourrainTrebuchet08,Stetter04}, logic and complexity \cite{CleggEdmondsImpagliazzo1996}, cryptography \cite{CourtoisKlimovPatarinShamir2000}, 
and optimization \cite{moniquefranz, Lasserre2002,Laurent2007, Parrilo2003, Parrilo2002, PokuttaSchulz2009}.

As the complexity of solving a combinatorial system with this strategy depends on its certificate degree, 
it is important to understand the class of problems having small degrees $D$.  In Theorem \ref{thm:non3color}, 
we give a combinatorial characterization of non-3-colorable graphs whose polynomial system encoding has a degree one Nullstellensatz certificate of infeasibility.
Essentially, a graph has a degree one certificate if there is an edge covering of
the graph by three and four cycles obeying some parity conditions on the number of times an edge is covered.
This result is reminiscent of the cycle double cover conjecture of Szekeres (1973)
\cite{Szekeres1973} and Seymour (1979) \cite{Seymour1979}.
The class of non-3-colorable graphs with degree one certificates is far from trivial;
it includes graphs that contain an odd-wheel or a 4-clique \cite{SusanPhd}
and experimentally it has been shown to include more complicated graphs
(see \cite{SusanPhd,DeLoeraLeeMalkinMargulies08,ima}).

%

In our second application of the polynomial method, we use tools from the theory of Gr\"obner bases to investigate (in Section \ref{hamiltonsec}) the detection of Hamiltonian cycles of a directed graph $G$.   The following ideals algebraically encode Hamiltonian cycles (see Lemma \ref{H_Gcountscycles} for a proof).

\begin{proposition} \label{openinghamiltonprop}
Let $G=(V,A)$ be a simple directed graph on vertices $V = \{1,\ldots,n\}$.
Assume that the characteristic of $\k$ is relatively prime to $n$ and that $\omega\in\k$ is a primitive $n$-th root of unity. Consider the following system in $\k[x_1,\ldots,x_n]$:
\[ H_G = \{  x_i^n-1=0,  \  \prod_{j \in \delta^+(i)}(\omega x_i - x_j)=0 :  \ i \in V\}.\]
Here, $\delta^+(i)$ denotes those vertices $j$ which are connected to $i$ by an arc going from $i$ to $j$ in $G$.  The system $H$ has a solution over $\overline{\k}$ if and only if $G$ has a Hamiltonian
cycle.
\end{proposition}

We prove a decomposition theorem for the ideal $H_G$ generated by the above polynomials, and 
based on this structure, we give an algebraic characterization of \emph{uniquely Hamiltonian graphs} 
(reminiscent of the one for $k$-colorability in \cite{HillarWindfeldt08}).  Our results also provide an algorithm to decide this property.  These findings are related to a well-known theorem of Smith \cite{Tutte1946} which states that if a $3$-regular graph has one Hamiltonian cycle then it has at least three.  
It is still an open question to decide the complexity of finding a second Hamiltonian cycle knowing that it exists \cite{Cameron2001}. 

Finally, in Section \ref{graphautosec} we explore the problem of determining the automorphisms $Aut(G)$ of an undirected graph $G$.   Recall that the elements of $Aut(G)$ are those permutations of the vertices of $G$ which preserve edge adjacency.  Of particular interest for us in that section is when graphs are \textit{rigid}; that is, $|Aut(G)| = 1$.  The complexity of this decision problem is still wide open \cite{cameron}. The combinatorial object $Aut(G)$ will be viewed as an algebraic variety in $\R^{n \times n}$ as follows. 

\begin{proposition}\label{AutGEncodingprop}
Let $G$ be a simple undirected graph and $A_G$ its adjacency matrix.  Then $Aut(G)$ is the group of permutation matrices $P = [P_{i,j}]_{i,j=1}^n$ given by the zeroes of the ideal $I_G \subseteq \R[x_1,\ldots,x_n]$ generated from the equations:
\begin{equation}\label{BinaryProgram1Intro}
\begin{split}
{(PA_G - A_GP)}_{i,j} = 0,& \ \ 1 \leq i,j \leq n; \ \ \  \sum_{i=1}^{n} P_{i,j} = 1, \ \ 1 \leq j \leq n; \\
\sum_{j=1}^{n} P_{i,j} = 1,& \  \ 1 \leq i \leq n; \ \ \ {P}^{2}_{i,j} - P_{i,j} = 0, \ \ 1 \leq i,j \leq n. \\
\end{split}
\end{equation}
\end{proposition}
\begin{proof}
The last three sets of equations say that $P$ is a permutation matrix, 
while the first one ensures that this permutation preserves adjacency of edges ($PA_GP^{\top} = A_G$).
\end{proof}

In what follows, we shall interchangeably refer to $Aut(G)$ as a group or the 
variety of Proposition \ref{AutGEncodingprop}. This real variety can be studied from the perspective of convexity. Indeed, 
from Proposition \ref{AutGEncodingprop}, $Aut(G)$ consists of the integer vertices of the polytope of doubly stochastic matrices commuting with $A_G$.  By replacing the equations ${P}^{2}_{i,j} - P_{i,j} = 0$ in (\ref{BinaryProgram1Intro}) with the linear inequalities $P_{ij} \geq 0$, we obtain a 
polyhedron $P_G$ which is a convex relaxation of the automorphism group of the graph. This polytope and its integer hull have been investigated by Friedland and Tinhofer \cite{tinhofer,friedlander}, where they gave conditions for it to be integral.  Here, we uncover more properties of the polyhedron $P_G$ and its integer vertices $Aut(G)$.

Our first result is that $P_G$ is \emph{quasi-integral}; that is, the graph induced by the integer points 
 in the 1-skeleton of $P_G$ is connected (see Definition 7.1 in Chapter 4 of \cite{Yemelichev}).  
It follows that one can decide rigidity of graphs by inspecting the vertex neighbors of the identity permutation.  
Another application of this result is an output-sensitive algorithm for enumerating all automorphisms of any graph \cite{avis}.
The problem of determining the triviality of the automorphism group of a graph can be solved efficiently when $P_G$ is integral. Such
graphs have been called \emph{compact} and a fair amount of research has been dedicated to them (see \cite{godsil,tinhofer} and references
therein). 
 
Next, we use the theory of Gouveia, Parrilo, and Thomas \cite{GouveiaParriloThomas2008}, applied to the ideal $I_G$ of Proposition \ref{AutGEncodingprop}, to approximate the integer hull of $P_G$ by projections of semidefinite programs
(the so-called \emph{theta bodies}). In their work, the authors of \cite{GouveiaParriloThomas2008} generalize the Lov\'asz theta body for $0/1$ polyhedra to generate a sequence of semidefinite programming relaxations computing the convex hull of the zeroes of a set of real polynomials \cite{LovaszSchrijver1991,Lovasz1994}.  The paper \cite{GouveiaParriloThomas2008} provides some applications to finding maximum stable sets \cite{LovaszSchrijver1991} and maximum cuts \cite{GouveiaParriloThomas2008}.   We study the theta bodies of the variety of automorphisms of a graph. 
In particular, we give sufficient conditions on $Aut(G)$ for which the first theta body is already equal to $P_G$ (in much the same way that stable sets of perfect 
graphs are theta-1 exact 
\cite{GouveiaParriloThomas2008, LovaszSchrijver1991}). Such graphs will be called \emph{exact}. Establishing these conditions for exactness requires an interesting generalization of properties of the symmetric group (see Theorem \ref{groebner} for details). 
In addition, we prove that compact graphs are a proper subset of exact graphs (see Theorem \ref{compactinsideexact}). 
This is interesting because we do not know of an example of a graph that is not exact, and the connection with semidefinite 
programming may open interesting approaches to understanding the complexity of the graph automorphism problem.

Below, we assume the reader is familiar with the basic properties of polynomial ideals 
and commutative algebra as introduced in the elementary text \cite{CoxLittleOShea92}. 
A quick, self-contained review can also be found in Section 2 of \cite{HillarWindfeldt08}.  

\section{Recognizing Non-3-colorable Graphs}\label{3coloringsec}

In this section, we give a complete combinatorial characterization of the class of non-3-colorable simple undirected graphs $G=(V,E)$ with a degree one Nullstellensatz certificate of infeasibility for the following system (with $\k = \F_2$) from Proposition \ref{prop:encodings}:
\begin{equation}
J_G = \{x_i^3+1=0, \ x_i^2+x_ix_j+x_j^2 = 0: \  i \in V, \ \{i,j\}\in E\}.
\label{3colorpolysys}
\end{equation}
This polynomial system has a degree one ($D = 1$) Nullstellensatz certificate of infeasibility if and only if there exist coefficients
$a_i,a_{ij},b_{ij},b_{ijk} \in \F_2$ such that
\begin{equation}
\sum_{i\in V} (a_i + \sum_{j\in V} a_{ij}x_j) (x_i^3+1)
+\sum_{\{i,j\}\in E} (b_{ij} + \sum_{k\in V} b_{ijk}x_k)(x_i^2+x_ix_j+x_j^2) 
= 1.
\label{deg1null}
\end{equation}

%

\begin{figure}[ht]
\begin{center}
\input{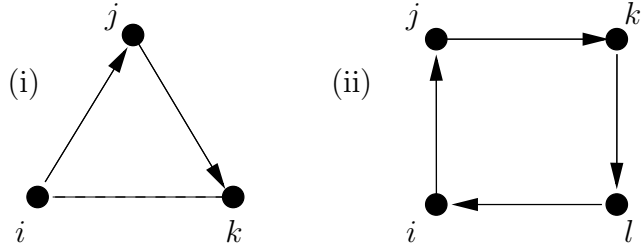}
\end{center}
\caption{(i) partial 3-cycle, (ii) chordless 4-cycle} 
\label{3and4cycles}
\end{figure}

Our characterization involves two types of substructures on the graph $G$ (see Figure \ref{3and4cycles}).  
The first of these are \emph{oriented partial-3-cycles}, which are pairs of arcs
$\{(i,j),(j,k)\} \subseteq Arcs(G)$, also denoted $(i,j,k)$, in which $(k,i)\in Arcs(G)$
(the vertices ${i,j,k}$ induce a 3-cycle in $G$).
The second are \emph{oriented chordless 4-cycles}, which are sets of four arcs
$\{(i,j),(j,k),(k,l),(l,i)\}  \subseteq Arcs(G)$, denoted $(i,j,k,l)$, with $(i,k),(j,l)\not\in Arcs(G)$
(the vertices ${i,j,k,l}$ induce a chordless 4-cycle).

\begin{theorem}
\label{thm:non3color}
For a given simple undirected graph $G=(V,E)$, the polynomial system
over $\F_2$ encoding the $3$-colorability of $G$
$$J_G = \{x_i^3+1=0, \ x_i^2+x_ix_j+x_j^2 = 0: \  i \in V, \ \{i,j\}\in E\}$$
has a degree one Nullstellensatz certificate of infeasibility if and
only if there exists a set $C$ of oriented partial $3$-cycles and oriented
chordless $4$-cycles from $Arcs(G)$ such that
\begin{enumerate}
\item $|C_{(i,j)}|+|C_{(j,i)}| \equiv 0 \pmod{2}$ for all $\{i,j\} \in E$ \ and \  \item $\sum_{(i,j)\in Arcs(G),i<j} |C_{(i,j)}| \equiv 1 \pmod{2}$,
\end{enumerate}
where $C_{(i,j)}$ denotes the set of cycles in $C$ in which the arc $(i,j) \in Arcs(G)$ appears.
Moreover, the class of non-3-colorable graphs whose encodings have
degree one Nullstellensatz infeasibility certificates can be recognized in
polynomial time.  
\end{theorem}


We can consider the set $C$ in Theorem~\ref{thm:non3color} as a covering of $E$ by directed edges.  From this perspective, Condition 1 in Theorem \ref{thm:non3color} means that every edge of $G$ is covered by an even number of arcs from cycles in $C$.  On the other hand, Condition 2 says that if $\hat{G}$ is the directed graph obtained from $G$ by the orientation induced by the total ordering on the vertices $1 < 2 < \cdots < n$, then when summing the number of times each arc in $\hat{G}$ appears in the cycles of $C$, the total is odd.  

Note that the 3-cycles and 4-cycles in $G$ that correspond to the partial
3-cycles and chordless 4-cycles in $C$ give an edge-covering 
of a non-3-colorable subgraph of $G$. 
Also, note that if a graph $G$ has a non-3-colorable subgraph whose polynomial
encoding has a degree one infeasibility certificate, then the encoding of
$G$ will also have a degree one infeasibility certificate.

The class of graphs with encodings that have degree one infeasibility certificates
includes all graphs containing odd wheels as subgraphs (e.g., a $4$-clique) \cite{SusanPhd}.
\begin{corollary}
If a graph $G=(V,E)$ contains an odd wheel, then the encoding of $3$-colorability of $G$ from Theorem \ref{thm:non3color} has
a degree one Nullstellensatz certificate of infeasibility.
\end{corollary}
\begin{proof}
Assume $G$ contains an odd wheel with vertices labelled as in Figure
\ref{fig:oddwheel} below.
Let $$C:=\{(i,1,i+1): 2\leq i\leq n-1\}\cup\{(n,1,2)\}.$$

\begin{figure}[ht]
\begin{center}
\includegraphics[scale=0.5]{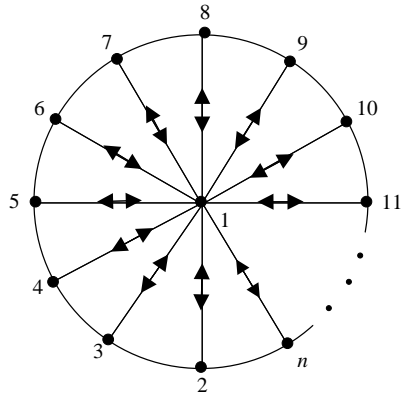}
\end{center}
\label{fig:oddwheel}
\caption{Odd wheel}
\end{figure}

Figure \ref{fig:oddwheel} illustrates the arc directions for the oriented partial
3-cycles of $C$.
Each edge of $G$ is covered by exactly zero or two partial 3-cycles, so
$C$ satisfies Condition 1 of Theorem \ref{thm:non3color}.
Furthermore, each arc $(1,i)\in Arcs(G)$ is covered exactly once by a partial
3-cycle in $C$, and there is an odd number of such arcs. Thus,
$C$ also satisfies Condition 2 of Theorem \ref{thm:non3color}.
\end{proof}

\begin{figure}[ht]
\begin{center}
\includegraphics[scale=0.2]{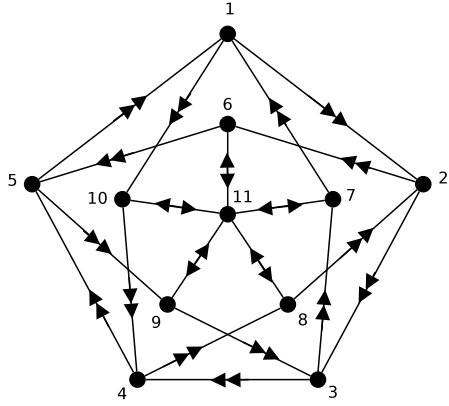}
\end{center}
\caption{Gr\"otzsch graph.} 
\label{groetzsch}
\end{figure}

A non-trivial example of a non-3-colorable graph with a degree one Nullstellensatz 
certicate is the Gr\"otzsch graph.
\begin{example}
Consider the Gr\"otzsch graph in Figure \ref{groetzsch}, which has no 3-cycles.  
The following set of oriented chordless 4-cycles gives a certificate of
non-3-colorability by Theorem~\ref{thm:non3color}: 
\begin{align*}
C:=\{ &(1,2,3,7), (2,3,4,8), (3,4,5,9),
(4,5,1,10), (1,10,11,7), \\ 
&(2,6,11,8),(3,7,11,9),(4,8,11,10),(5,9,11,6)\}.
\end{align*}
Figure \ref{groetzsch} illustrates the arc directions for the 4-cycles of $C$.
Each edge of the graph is covered by exactly two 4-cycles, so $C$ satisfies
Condition 1 of Theorem \ref{thm:non3color}.
Moreover, one can check that Condition 2 is also satisfied. It follows that the graph has no proper 3-coloring.
\qed
\end{example}

We now prove Theorem~\ref{thm:non3color} using ideas from polynomial algebra.
First, notice that we can simplify a degree one certificate as follows:
Expanding the left-hand side of (\ref{deg1null}) and collecting terms, the
only coefficient of $x_jx_i^3$ is $a_{ij}$ and thus $a_{ij}=0$ for all $i,j\in
V$. Similarly, the only coefficient of $x_ix_j$ is $b_{ij}$, and so
$b_{ij}=0$ for all $\{i,j\}\in E$.  We thus arrive at the following simplified
expression:
\begin{equation} \sum_{i\in V} a_i(x_i^3+1)
+\sum_{\{i,j\}\in E} (\sum_{k\in V} b_{ijk}x_k)(x_i^2+x_ix_j+x_j^2) 
= 1.\end{equation}
Now, consider the following set $F$ of polynomials:
\begin{align}
& x_i^3+1 & \qquad  \forall i\in V, \\
& x_k(x_i^2+x_ix_j+x_j^2) & \qquad \forall \{i,j\}\in E, \ k \in V. 
\label{eqn:edge}
\end{align}
The elements of $F$ are those polynomials that can appear in a degree one certificate of infeasibility. Thus, there exists a degree one certificate if and only if the constant
polynomial 1 is in the linear span of $F$; that is, $1\in\gen{F}$, where
$\gen{F}$ is the vector space over $\F_2$ generated by the polynomials in $F$.

We next simplify the set $F$. Let $H$ be the following set of polynomials:
\begin{align}
&x_i^2x_j+x_ix_j^2+1 &\; \forall \{i,j\}\in E, \label{eqn:Hfirstpolys} \\
&x_ix_j^2+x_jx_k^2 &\;  \forall (i,j),(j,k),(k,i)\in Arcs(G),
\label{eqn:3cycle}\\
&x_ix_j^2+x_jx_k^2+x_kx_l^2+x_lx_i^2 &\;  \forall
(i,j),(j,k),(k,l),(l,i)\in Arcs(G), (i,k),(j,l)\not\in Arcs(G). \label{eqn:4cycle}
\end{align}
If we identify the monomials $x_ix_j^2$ as the arcs $(i,j)$, then the polynomials (\ref{eqn:3cycle}) correspond to oriented partial 3-cycles and the polynomials (\ref{eqn:4cycle}) correspond to oriented chordless 4-cycles. The following lemma says that we can use $H$ instead of $F$ to find a degree one certificate.
\begin{lemma}
\label{lem:deg1equiv}
We have $1 \in \gen{F}$ if and only if $1 \in  \gen{H}$.
\end{lemma}
\begin{proof}
The polynomials (\ref{eqn:edge}) above can be split into two classes of
equations:
(i) $k=i$ or $k=j$ and (ii) $k\ne i$ and $k\ne j$.
Thus, the set $F$ consists of
\begin{align}
& x_i^3+1 \qquad & \forall i\in V, \label{eqn:vertex} \\
& x_i(x_i^2+x_ix_j+x_j^2)=x_i^3+x_i^2x_j+x_ix_j^2 & \qquad  \forall \{i,j\}\in E,
\label{eqn:edge0} \\
& x_k(x_i^2+x_ix_j+x_j^2)=x_i^2x_k+x_ix_jx_k+x_j^2x_k & \qquad
\forall  \{i,j\}\in E, \ k \in V, i\ne k \ne j.
\end{align}
Using polynomials (\ref{eqn:vertex}) to eliminate the $x_i^3$ terms from (\ref{eqn:edge0}), we arrive at the following set of polynomials, which we label $F'$:
\begin{align}
& x_i^3+1 \qquad &  \forall i\in V, \label{eqn:vertex2}\\
& x_i^2x_j+x_ix_j^2+1 = (x_i^3+x_i^2x_j+x_ix_j^2)+(x_i^3+1)
& \qquad  \forall \{i,j\}\in E, \label{eqn:edge1}\\
& x_i^2x_k+x_ix_jx_k+x_j^2x_k & \qquad \forall \{i,j\}\in E, \ k \in V, i\ne k \ne j.
\end{align}
Observe that $\gen{F}=\gen{F'}$. We can eliminate the polynomials
(\ref{eqn:vertex2}) as follows.  For every 
$i\in V$, $(x_i^3+1)$ is the only polynomial in $F'$ containing the monomial
$x_i^3$ and thus the polynomial $(x_i^3+1)$ cannot be present in any nonzero 
linear combination of the polynomials in $F'$ that equals 1.
We arrive at the following smaller set of polynomials, which we label $F''$.
\begin{align}
& x_i^2x_j+x_ix_j^2+1 & \qquad  \forall \{i,j\}\in E, \label{eqn:edge3}\\
& x_i^2x_k+x_ix_jx_k+x_j^2x_k & \qquad  \forall \{i,j\}\in E, k \in V, i\ne k \ne j.
\end{align}
So far, we have shown $1\in\gen{F}=\gen{F'}$ if and only if $1\in\gen{F''}$.

Next, we eliminate monomials of the form $x_ix_jx_k$. There are 3 cases to consider.

Case 1: $\{i,j\}\in E$ but $\{i,k\}\not\in E$ and $\{j,k\}\not\in E$.
In this case, the monomial $x_ix_jx_k$ appears in only one polynomial,
$x_k(x_i^2+x_ix_j+x_j^2)=x_i^2x_k+x_ix_jx_k+x_j^2x_k$, so we can eliminate all
such polynomials.

Case 2: $i,j,k\in V$, $(i,j),(j,k),(k,i) \in Arcs(G)$.
Graphically, this represents a 3-cycle in the graph. In this case, the monomial $x_ix_jx_k$ appears in three polynomials:
\begin{align}
& x_k(x_i^2+x_ix_j+x_j^2)=x_i^2x_k+x_ix_jx_k+x_j^2x_k, \label{eqn:case21}\\
& x_j(x_i^2+x_ix_k+x_k^2)=x_i^2x_j+x_ix_jx_k+x_jx_k^2, \\
& x_i(x_j^2+x_jx_k+x_k^2)=x_ix_j^2+x_ix_jx_k+x_ix_k^2.
\end{align}
Using the first polynomial, we can eliminate $x_ix_jx_k$ from the other two:
\begin{align*}
& x_i^2x_j+ x_jx_k^2 + x_i^2x_k +x_j^2x_k
= (x_i^2x_j+x_ix_jx_k+x_jx_k^2) + (x_i^2x_k+x_ix_jx_k+x_j^2x_k),\\
& x_ix_j^2 + x_ix_k^2 + x_i^2x_k + x_j^2x_k
= (x_ix_j^2+x_ix_jx_k+x_ix_k^2) +  (x_i^2x_k+x_ix_jx_k+x_j^2x_k).
\end{align*}
We can now eliminate the polynomial (\ref{eqn:case21}).
Moreover, we can use the polynomials (\ref{eqn:edge3}) 
to rewrite the above two polynomials as follows.
\begin{align*}
&x_kx_i^2+x_ix_j^2= (x_i^2x_j+x_jx_k^2+x_i^2x_k+x_j^2x_k)
+(x_jx_k^2+x_j^2x_k+1) + (x_ix_j^2+x_i^2x_j+1), \\
& x_ix_j^2+x_jx_k^2 = (x_ix_j^2+x_ix_k^2+x_i^2x_k+x_j^2x_k) 
+(x_ix_k^2+x_i^2x_k+1)+(x_jx_k^2+x_j^2x_k+1).
\end{align*}
Note that both of these polynomials correspond to two of the arcs of the $3$-cycle 
$(i,j),(j,k),(k,i) \in Arcs(G)$.

Case 3: $i,j,k\in V$, $(i,j),(j,k)\in Arcs(G)$ and $(k,i)\not\in Arcs(G)$.  We have
\begin{align}
& x_k(x_i^2+x_ix_j+x_j^2)=x_i^2x_k+x_ix_jx_k+x_j^2x_k \label{eqn:case31}, \\
& x_i(x_j^2+x_jx_k+x_k^2)=x_ix_j^2+x_ix_jx_k+x_ix_k^2.
\end{align}
As before we use the first polynomial to eliminate the monomial $x_ix_jx_k$
from the second:
\begin{align*}
x_ix_j^2+x_jx_k^2+(x^2_ix_k+x_ix_k^2+1) = \ &
(x_ix_j^2+x_ix_jx_k+x_ix_k^2)+ (x_i^2x_k+x_ix_jx_k+x_j^2x_k) \\
&+ (x_jx_k^2+x_j^2x_k+1).
\end{align*}
We can now eliminate (\ref{eqn:case31}); thus, the original system has been reduced to the following one, which we label as $F'''$:
\begin{align}
&x_i^2x_j+x_ix_j^2+1 &\;  \forall \{i,j\}\in E, \\
&x_ix_j^2+x_jx_k^2  &\;   \forall (i,j),(i,k),(j,k)\in Arcs(G),\\
&x_ix_j^2+x_jx_k^2+(x^2_ix_k+x_ix_k^2+1) &\;   \forall (i,j),(j,k)\in Arcs(G), (k,i)\not\in Arcs(G). \label{eqn:v}
\end{align}
Note that $1\in\gen{F}$ if and only if $1\in\gen{F'''}$.

The monomials $x^2_ix_k$ and $x_ix_k^2$ with $(k,i)\not\in Arcs(G)$ always appear
together and only in the polynomials (\ref{eqn:v}) in the expression $(x^2_ix_k+x_ix_k^2+1)$. Thus, we can eliminate  
the monomials $x^2_ix_k$ and $x_ix_k^2$ with $(k,i)\not\in Arcs(G)$ by choosing one of the polynomials (\ref{eqn:v}) and using it to eliminate
the expression $(x^2_ix_k+x_ix_k^2+1)$ from all other polynomials in which it appears.
Let $i,j,k,l\in V$ be such that $(i,j),(j,k),(k,l),(l,i)\in Arcs(G)$ and $(k,i),(i,k)\not\in Arcs(G)$.
We can then eliminate the monomials $x^2_ix_k$ and $x_ix_k^2$ as follows:
\begin{align*} 
x_ix_j^2+x_jx_k^2+x_kx_l^2+x_lx_i^2
= & \ (x_ix_j^2+x_jx_k^2+x^2_ix_k+x_ix_k^2+1) \\&+ (x_kx_l^2+x_lx_i^2+x^2_ix_k+x_ix_k^2+1).
\end{align*}
Finally, after eliminating the polynomials (\ref{eqn:v}), we have system $H$ (polynomials (\ref{eqn:Hfirstpolys}),
(\ref{eqn:3cycle}), and (\ref{eqn:4cycle})):
\begin{align*}
&x_i^2x_j+x_ix_j^2+1 &\; \forall \{i,j\}\in E, \\
&x_ix_j^2+x_jx_k^2 &\;  \forall (i,j),(j,k),(k,i)\in Arcs(G),
\label{eqn:3cycle}\\
&x_ix_j^2+x_jx_k^2+x_kx_l^2+x_lx_i^2 &\;  \forall
(i,j),(j,k),(k,l),(l,i)\in Arcs(G), (i,k),(j,l)\not\in Arcs(G).
\end{align*}
The system $H$ has the property that $1\in\gen{F'''}$ if and only if $1\in\gen{H}$, and thus, $1\in \gen{F}$ if and only if $1\in\gen{H}$ as required
\end{proof}

We now establish that the sufficient condition for infeasibility $1\in\gen{H}$ is equivalent to
the combinatorial parity conditions in Theorem \ref{thm:non3color}.

%

\begin{lemma}
\label{lem:deg1equiv2}
There exists a set $C$ of oriented partial 3-cycles and oriented chordless
4-cycles satisfying Conditions 1. and 2. of Theorem~\ref{thm:non3color} if and only if
$1\in\gen{H}$.
\end{lemma}
\begin{proof}
Assume that $1\in\gen{H}$.
Then there exist coefficients $c_h\in\F_2$ such that $\sum_{h\in H} c_hh = 1$.
Let $H':=\{h\in H: c_h =1 \}$; then, $\sum_{h\in H'}h = 1$.
Let $C$ be the set of oriented partial 3-cycles $(i,j,k)$ where
$x_ix_j^2+x_jx_k^2 \in H'$ together with the set of oriented chordless
4-cycles $(i,j,l,k)$ where 
$x_ix_j^2+x_jx_l^2+x_lx_k^2+x_kx_i^2 \in H'$.
Now, $|C_{(i,j)}|$ is the number of polynomials in $H'$ of the form
(\ref{eqn:3cycle}) or (\ref{eqn:4cycle}) in which the monomial $x_ix_j^2$ appears,
and similarly, $|C_{(j,i)}|$ is the number of polynomials in $H'$ of the form
(\ref{eqn:3cycle}) or (\ref{eqn:4cycle}) in which the monomial $x_jx_i^2$ appears.
Thus, $\sum_{h\in H'} h = 1$ implies that, for every pair $x_ix_j^2$
and $x_jx_i^2$, either
\begin{enumerate}
\item $|C_{(i,j)}| \equiv 0 \pmod{2}$, $|C_{(j,i)}| \equiv 0 \pmod{2}$,
and $x_i^2x_j+x_ix_j^2+1\not\in H'$ or
\item $|C_{(i,j)}| \equiv 1 \pmod{2}$, $|C_{(j,i)}| \equiv 1 \pmod{2}$,
and $x_i^2x_j+x_ix_j^2+1\in H'$.
\end{enumerate}
In either case, we have $|C_{(i,j)}|+|C_{(j,i)}|\equiv 0\pmod{2}$.
Moreover, since $\sum_{h\in H'} h = 1$, there must be an odd number of the
polynomials of the form $x_i^2x_j+x_ix_j^2+1$ in $H'$.  That is, case 2 above 
occurs an odd number of times and therefore, $\sum_{(i,j)\in Arcs(G),i<j} |C_{(i,j)}| \equiv 1 \pmod{2}$ as required.

Conversely, assume that there exists a set $C$ of oriented partial 3-cycles and
oriented chordless 4-cycles satisfying the conditions of Theorem~\ref{thm:non3color}.
Let $H'$ be the set of polynomials $x_ix_j^2+x_jx_k^2$ 
where $(i,j,k)\in C$ and the set of polynomials
$x_ix_j^2+x_jx_l^2+x_lx_k^2+x_kx_i^2$ where 
$(i,j,l,k)\in C$ together with 
the set of polynomials $x_i^2x_j+x_ix_j^2+1\in H$ where $|C_{(i,j)}| \equiv 1$.
Then, $|C_{(i,j)}|+|C_{(j,i)}|\equiv 0\pmod{2}$ implies that
every monomial $x_ix_j^2$ appears in an even number polynomials of $H'$. 
Moreover, since $\sum_{(i,j)\in Arcs(G),i<j} |C_{(i,j)}| \equiv 1 \pmod{2}$, there 
are an odd number of polynomials $x_i^2x_j+x_ix_j^2+1$ appearing in $H'$.
Hence, $\sum_{h\in H'} h = 1$ and $1\in\gen{H}$.
\end{proof}

Combining Lemmas \ref{lem:deg1equiv} and \ref{lem:deg1equiv2},
we arrive at the characterization stated in Theorem~\ref{thm:non3color}.
That such graphs can be decided in polynomial time follows from the fact that the existence of a certificate of any fixed degree  can be decided in polynomial time (as is well known and follows since there are polynomially many monomials up to any fixed degree; see also \cite[Theorem 4.1.3]{SusanPhd}).

Finally, we pose as open problems the construction of a variant of Theorem \ref{thm:non3color} for general $k$-colorability and also combinatorial characterizations for larger certificate degrees $D$. 

\begin{problem}
Characterize those graphs with a given $k$-colorability Nullstellensatz certificate of degree $D$.   
\end{problem}

\section{Recognizing Uniquely Hamiltonian Graphs} \label{hamiltonsec}

Throughout this section we work over an arbitrary algebraically 
closed field $\k = \kc$, although in some cases, we will need to restrict its characteristic. 
Let us denote by $H_G$ the \textit{Hamiltonian ideal} generated by the polynomials from Proposition \ref{openinghamiltonprop}.  A connected, directed graph $G$ with $n$ vertices has a Hamiltonian
cycle if and only if the equations defined by $H_{G}$ have a solution over
$\k$ (or, in other words, if and only if $V(H_{G}) \neq \emptyset$ for the algebraic variety $V(H_{G})$ associated to the ideal $H_G$).  In a precise sense to be made clear below, the ideal $H_G$ actually encodes \textit{all} Hamiltonian cycles of $G$.  
However, we need to be somewhat careful about how to count cycles (see Lemma \ref{H_Gcountscycles}). 
In practice  $\omega$ can be treated as a variable and not as a fixed primitive $n$-th root of unity. 
A set of equations ensuring that $\omega$ only takes on the value of a \emph{primitive} $n$-th root of unity is the following: \[\{\omega^{i(n-1)} + \omega^{i(n-2)} + \cdots + \omega^i + 1 = 0:  \ 1 \leq i \leq n\}.\]
We can also use the cyclotomic polynomial $\Phi_n(\omega)$ \cite{dummit}, which is the polynomial whose zeroes are the
primitive $n$-th roots of unity. 

%

We shall utilize the theory of Gr\"obner bases to show that $H_G$ has a special (algebraic) decomposition structure in terms of the different Hamiltonian cycles of $G$ (this is Theorem \ref{H_Gdecompthm} below).  In the particular case when $G$ has a unique Hamiltonian cycle, we get a specific algebraic criterion which can be algorithmically verified.  These results are Hamiltonian analogues to the algebraic $k$-colorability characterizations of \cite{HillarWindfeldt08}.  We first turn our attention more generally to cycle ideals of a simple directed graph $G$.  These will be the basic elements in our decomposition of the Hamiltonian ideal $H_G$, as they algebraically encode single cycles $C$ (up to symmetry).

When $G$ has the property that each pair of vertices connected by an arc is also connected by an arc in the opposite direction, then we call $G$ \emph{doubly covered}.  When $G = (V,E)$ is presented as an undirected graph, we shall always view it as the doubly covered directed graph on vertices $V$ with arcs $Arcs(G)$.

Let $C$ be a cycle of length $k > 2$ in $G$, expressed as a sequence of arcs, \[C = \{ (v_1,v_2),(v_2,v_3),\ldots,(v_k,v_1)\}.\]  For the purpose of this work, we call $C$ a \emph{doubly covered cycle} if consecutive vertices in the cycle are connected by arcs in both directions; otherwise, $C$ is simply called \emph{directed}.  In particular, each cycle in a doubly covered graph is a doubly covered cycle.  These definitions allow us to work with both undirected and directed graphs in the same framework.

\begin{definition}[Cycle encodings]
Let $\omega$ be a fixed primitive $k$-th root of unity and let $\k$ be a field with characteristic not dividing $k$.  If $C$ is a doubly covered cycle of length $k$ and the vertices in $C$ are $\{v_1,\ldots,v_k\}$, then the \textit{cycle encoding} of $C$ is the following set of $k$ polynomials
in $\k[x_{v_1},\ldots,x_{v_k}]$:
\begin{equation}\label{cycleGBundirected}
g_i  = \begin{cases} 
\  x_{v_i} + \frac{(\omega^{2+i} - \omega^{2-i})}{(\omega^3-\omega)} x_{v_{k-1}} + \frac{(\omega^{1-i} - \omega^{3+i})}{(\omega^3-\omega)}x_{v_k}   &  \ \ \text{$i = 1,\ldots,k-2$}, \\
\  (x_{v_{k-1}} - \omega x_{v_k})(x_{v_{k-1}} - \omega^{-1} x_{v_k})   & \ \ \text{$i = k-1$},  \\
\  x_{v_k}^k -1  & \ \ \text{$i = k$}.
 \end{cases}
\end{equation}

If $C$ is a directed cycle of length $k$ in a directed graph, with vertex set  $\{v_1,\ldots,v_k\}$, the \textit{cycle encoding} of $C$ is the following set of $k$ polynomials:
\begin{equation}\label{cycleGBdirected}
g_i  = \begin{cases} 
\  x_{v_{k-i}} - \omega^{k-i} x_{v_k}  &  \ \ \text{$i = 1,\ldots,k-1$}, \\
\  x_{v_k}^k -1  & \ \ \text{$i = k$}.
 \end{cases}
\end{equation}
\end{definition}

\begin{definition}[Cycle Ideals]
The \textit{cycle ideal} associated to a cycle $C$ is
$H_{G,C} = \<g_1,\ldots,g_k \> \subseteq \k[x_{v_1},\ldots,x_{v_k}]$,
where the $g_i$s are the cycle encoding of $C$ given by (\ref{cycleGBundirected}) or (\ref{cycleGBdirected}).
\end{definition}

%

The polynomials $g_i$ are computationally useful generators for cycle ideals.
(Once again, see \cite{CoxLittleOShea92} for the relevant background on
Gr\"obner bases and term orders.)

\begin{lemma}\label{cycleGB}
The cycle encoding polynomials $F = \{g_1,\ldots,g_k\}$ are a reduced Gr\"obner basis for the cycle ideal $H_{G,C}$ with respect to any term order $\prec$ with $x_{v_k} \prec \cdots  \prec x_{v_1}$.
\end{lemma}
\begin{proof}
Since the leading monomials in a cycle encoding:
\begin{equation}\label{leadingmonomials}
\{x_{v_1},\ldots,x_{v_{k-2}},x_{v_{k-1}}^2,x_{v_k}^k\} \ \textit{ \rm or } \ \{ x_{v_1},\ldots,x_{v_{k-2}},x_{v_{k-1}},x_{v_k}^k \}
\end{equation}
are relatively prime, the polynomials $g_i$ form a Gr\"obner
basis for $H_{G,C}$ (see Theorem 3 and Proposition 4 in \cite[Section 2]{CoxLittleOShea92}).  That $F$ is reduced follows from inspection of (\ref{cycleGBundirected}) and (\ref{cycleGBdirected}).  
\end{proof}
\begin{remark}
In particular, since reduced Gr\"obner bases (with respect to a fixed term order) are unique, it follows that cycle encodings are canonical ways of generating cycle ideals (and thus of representing cycles by Lemma \ref{cycencodingthm}). 
\end{remark}

Having explicit Gr\"obner bases for these ideals allows us to compute their Hilbert series easily.

\begin{corollary}
The Hilbert series of $\k[x_{v_1},\ldots,x_{v_k}]/H_{G,C}$ for a doubly covered cycle or a directed cycle is equal to (respectively) \[ \frac{(1-t^2)(1-t^k)}{(1-t)^2} \text{ or } \ \frac{(1-t^k)}{(1-t)}. \]
\end{corollary}
\begin{proof}
If $\prec$ is a graded term order, then the (affine) Hilbert function of an ideal and of its ideal of leading terms are the same \cite[Chapter 9, \S 3]{CoxLittleOShea92}.  The form of the Hilbert series is now immediate from (\ref{leadingmonomials}).
\end{proof}

The naming of these ideals is motivated by the following result; in words, it 
says that the cycle $C$ is encoded as a complete intersection by the ideal $H_{G,C}$.

\begin{lemma}\label{cycencodingthm}
The following hold for the ideal $H_{G,C}$.
\begin{enumerate}
\item $H_{G,C}$ is radical, 
\item $|V(H_{G,C})| = k$ if $C$ is directed, and $|V(H_{G,C})| = 2k$ if $C$ is doubly covered undirected.
\end{enumerate}
\end{lemma}

\begin{proof}
Without loss of generality, we suppose that $v_i = i$ for $i = 1,\ldots,k$.  Let $\prec$ 
be any term order in which $x_k \prec \cdots  \prec x_1$.  From Lemma \ref{cycleGB}, the set of $g_i$ form a Gr\"obner basis for $H_{G,C}$.  It follows that the number of standard monomials of $H_{G,C}$ 
is $2k$ if $C$ is doubly covered undirected (resp. $k$ if it is directed).  Therefore by \cite[Lemma 2.1]{HillarWindfeldt08}, if we can prove that $|V(H_{G,C})| \geq k$ (resp. $|V(H_{G,C})| \geq 2k$), then both statements 1. and 2. follow.  

When $C$ is directed, this follows easily from the form of (\ref{cycleGBdirected}), so we shall assume that $C$ is doubly covered undirected.  We claim that the $k$ cyclic permutations of the two points:
\[ (\omega,\omega^2,\ldots,\omega^k), (\omega^k,\omega^{k-1},\ldots,\omega)\]
are zeroes of $g_i$, $i = 1,\ldots,k$.  Since cyclic permutation is multiplication by
a power of $\omega$, it is clear that we need only verify this claim for the two points above.
In the fist case, when $x_i = \omega^i$, we compute that for $i = 1,\ldots,k-2$:
\begin{equation*}
\begin{split}
(\omega^3-\omega) g_i(\omega,\ldots,\omega^k) = \ & (\omega^3-\omega) \omega^i + (\omega^{2+i} - 
\omega^{2-i}) \omega^{k-1} + (\omega^{1-i} - \omega^{3+i}) \omega^{k} \\
= \ & \omega^{3+i}-\omega^{1+i} + \omega^{1+i+k} - 
\omega^{1-i+k} + \omega^{1-i+k} - \omega^{3+i+k}  \\
= \ & 0,
\end{split}
\end{equation*}
since $\omega^k = 1$.  In the second case, when $x_i = \omega^{1-i}$, we 
again compute that for $i = 1,\ldots,k-2$:
\begin{equation*}
\begin{split}
(\omega^3-\omega) g_i(\omega^k,\ldots,\omega) = \ & (\omega^3-\omega) \omega^{1-i} + (\omega^{2+i} - 
\omega^{2-i}) \omega^{2} + (\omega^{1-i} - \omega^{3+i}) \omega \\
= \ & \omega^{4-i}-\omega^{2-i} + \omega^{4+i} - 
\omega^{4-i} + \omega^{2-i} - \omega^{4+i}  \\
= \ & 0.
\end{split}
\end{equation*}
Finally, it is obvious that the two points zero $g_{k-1}$ and $g_k$, and this completes
the proof.
\end{proof}
\begin{remark}
Conversely, it is easy to see that points in $V(H_{G,C})$ correspond to cycles of length $k$ in $G$.  That this variety contains $k$ or $2k$ points corresponds to there being $k$ or $2k$ ways of writing down the cycle since we may cyclically permute it and also reverse its orientation (if each arc in the path is bidirectional).
\end{remark}

Before stating our decomposition theorem (Theorem \ref{H_Gdecompthm}), we need to explain how the Hamiltonian ideal encodes all Hamiltonian cycles of the graph $G$.  

\begin{lemma}\label{H_Gcountscycles}
Let $G$ be a connected directed graph on $n$ vertices.  Then, 
\[ V(H_G) = \bigcup_C V(H_{G,C}),\]
where the union is over all Hamiltonian cycles $C$ in $G$.  
\end{lemma}
\begin{proof}
We only need to verify that points in $V(H_G)$ correspond to cycles of length $n$.  Suppose there exists a Hamiltonian cycle in the graph $G$.  Label vertex $1$ in the cycle with the number $x_1 = \omega^0 = 1$ and then successively label vertices along the cycle with one higher power of $\omega$.  It is clear that these labels $x_i$ associated to vertices $i$ zero all of the equations generating $H_G$.

Conversely, let $\textbf{v} = (x_1,\ldots,x_n)$ be a point in the variety $V(H_G)$ associated to $H_G$; we claim that $\textbf{v}$ encodes a Hamiltonian cycle.  From the edge equations, each vertex must be adjacent to one labeled with the next highest power of $\omega$.  Fixing a starting vertex $i$, it follows that there is a cycle $C$ labeled with (consecutively) increasing powers of $\omega$.  Since $\omega$ is a primitive $n$th root of unity, this cycle must have length $n$, and thus is Hamiltonian.  
\end{proof}

Combining all of these ideas, we can prove the following result.

\begin{theorem}\label{H_Gdecompthm}
Let $G$ be a connected directed graph with $n$ vertices.  Then,
\[ H_G = \bigcap_{C} H_{G,C},\]
where $C$ ranges over all Hamiltonian cycles of the graph $G$.
\end{theorem}
\begin{proof}
Since $H_G$ contains a square-free univariate polynomial in each indeterminate, it is radical (see for instance \cite[Lemma 2.1]{HillarWindfeldt08}).  It follows that 
\begin{equation}
\begin{split}
 H_G = \ & I(V(H_G))  \\
= \ & I \left( \bigcup_C V(H_{G,C}) \right) \\
= \ & \bigcap_C I(V(H_{G,C}) )\\
= \ & \bigcap_C H_{G,C},\\
\end{split}
\end{equation}
where the second inequality comes from Lemma \ref{H_Gcountscycles} and the last one from  $H_{G,C}$ being a radical ideal (Lemma \ref{cycencodingthm}).
\end{proof}

We call a directed graph (resp. doubly covered graph) \textit{uniquely Hamiltonian} if it contains $n$ cycles of length $n$ (resp. $2n$ cycles of length $n$).

\begin{corollary} The graph $G$ is uniquely Hamiltonian if and only if the Hamiltonian ideal $H_G$ is of the form  $H_{G,C}$ for some length $n$ cycle $C$.
\end{corollary}

This corollary provides an algorithm to check whether a graph is uniquely Hamiltonian.
We simply compute a unique reduced Gr\"{o}bner basis of $H_G$ and then check that it has the same form as that of an ideal $H_{G,C}$.  Another approach is to count the number of standard monomials of any Gr\"obner bases for $H_G$ and compare with $n$ or $2n$ (since $H_G$ is radical).  We remark, however, that it is well-known that computing a Gr\"obner basis in general cannot be done in polynomial time \cite[p. 400]{Yap2000}.

We close this section with a directed and an undirected  example of Theorem \ref{H_Gdecompthm}.

\begin{example}
Let $G$ be the directed graph with vertex set $V = \{ 1,2,3,4,5 \}$ and arcs $A = \{ (1,2),(2,3),(3,4),(4,5),(5,1),(1,3),(1,4) \}$.  Moreover, let $\omega$ be a primitive $5$-th root of unity.  The ideal $H_G \subset \k[x_1,x_2,x_3,x_4,x_5]$ is generated by the polynomials,
\[  \{ x_i^5 - 1 : 1 \leq i \leq 5\} \cup \{(\omega x_1 - x_2)(\omega x_1 - x_3)(\omega x_1 - x_4), \omega x_2 - x_3, \omega x_3 - x_4, \omega x_4 - x_5, \omega x_5 - x_1\}. \]
A reduced Gr\"{o}bner basis for $H_G$ with respect to the ordering $x_5 \prec x_4 \prec x_3 \prec x_2 \prec x_1$ is
\[ \{ x_5^5-1, x_4 - \omega^4 x_5, x_3 - \omega^3 x_5, x_2 - \omega^2 x_5, x_1 - \omega x_5 \},\]
which is a generating set for $H_{G,C}$ with $C = \{(1,2), (2,3), (3,4),(4,5), (5,1)\}$.  \qed
\end{example}

Let $G$ be an undirected graph with vertex set $V$ and edge set $E$, and consider the auxiliary directed graph $\tilde{G}$ with vertices $V$ and arcs $Arcs(G)$.  Notice that $\tilde{G}$ is doubly covered, and hence each of its cycles are doubly covered.  We apply Theorem~\ref{H_Gdecompthm} to $H_{ \tilde{G} }$ to determine and count Hamiltonian cycles in $G$.  In particular, the cycle $C = \{v_1,v_2,\ldots,v_n\}$ of $G$ is Hamiltonian if and only if $\{(v_1, v_2), (v_2, v_3), \ldots, (v_{n-1}, v_n), (v_n, v_1) \}$ and $\{(v_2, v_1), (v_3, v_2), \ldots, (v_n, v_{n-1}), (v_1, v_n) \}$ are Hamiltonian cycles of $\tilde{G}$.

\begin{example}
Let $G$ be the undirected complete graph on the vertex set $V=\{1,2,3,4\}$.  Let $\tilde{G}$ be the doubly covered graph with vertex set $V$ and arcs $Arcs(G)$. 
Notice that $\tilde{G}$ has twelve Hamiltonian cycles:
\begin{align*}
C_1 =& \{(1,2),(2,3),(3,4),(4,1)\},& \ \ C_2 =& \{(2,1),(3,2),(4,3),(1,4)\}, &\\
C_3 =& \{(1,2),(2,4),(4,3),(3,1)\},& \ \ C_4 =& \{(2,1),(4,2),(3,4),(1,3)\}, &\\
C_5 =& \{(1,3),(3,2),(2,4),(4,1)\},& \ \ C_6 =& \{(3,1),(2,3),(4,2),(1,4)\}, &\\
C_7 =& \{(1,3),(3,4),(4,2),(2,1)\},& \ \ C_8 =& \{(3,1),(4,3),(2,4),(1,2)\}, &\\
C_9 =& \{(1,4),(4,2),(2,3),(3,1)\},& \ \ C_{10} =& \{(4,1),(2,4),(3,2),(1,3)\}, &\\
C_{11} =& \{(1,4),(4,3),(3,2),(2,1)\},& \ \ C_{12} =& \{(4,1),(3,4),(2,3),(1,2)\}. 
\end{align*}
One can check in a symbolic algebra system such as SINGULAR or Macaulay 2 that the ideal $H_{\tilde{G}}$ is the intersection
of the cycle ideals $H_{\tilde{G},C_i}$ for $i = 1,\ldots,12$.
\end{example}

\section{Permutation Groups as Algebraic Varieties and their Convex Approximations}
\label{graphautosec}

In this section, we study convex hulls of permutations groups viewed as permutation matrices.
We begin by studying the convex hull of automorphism groups of undirected simple graphs; these have 
a natural polynomial presentation using Proposition~\ref{AutGEncodingprop} from the introduction.
For background material on graph automorphism groups see \cite{cameron,godsil}.

We write $Aut(G)$ for the automorphism group of a graph $G = (V,E)$.  Elements of $Aut(G)$ are naturally represented as $|V| \times |V|$ 
permutation matrices; they are the \emph{integer} vertices of the rational polytope $P_G$ defined in the discussion following Proposition~\ref{AutGEncodingprop}. The polytope $P_G$ was first introduced by Tinhofer \cite{tinhofer}.   
Since we are primarily interested in the integer vertices of $P_G$, we investigate $IP_G$, the \textit{integer hull} of $P_G$ (i.e. $IP_G = conv(P_G \cap \Z^{n \times n})$). In the fortunate case that $P_G$ is already integral ($P_G$ = $IP_G$), we say that the graph $G$ is \emph{compact}, a term coined in \cite{tinhofer}.  This occurs, for example, in the special case that $G$ is an independent set on $n$ vertices. In this case $Aut(G)=S_n$ and  $P_G$ is the well-studied Birkhoff polytope, the convex hull of all doubly-stochastic matrices (see Chapter 5 of \cite{Yemelichev}).  One can therefore view $P_G$ as a generalization of the Birkhoff polytope to arbitrary graphs. Unfortunately, the polytope $P_G$ is not always integral. For instance, $P_G$ is not integral when $G$ is the Petersen graph.  Nevertheless, we can prove the following related result.

\begin{proposition}
The polytope $P_G$ is quasi-integral.  That is, the induced subgraph of the integer points of the 1-skeleton of $P_G$ is connected.
\end{proposition}

\begin{proof}
We claim that there exists a $0/1$ matrix $A$ such that $P_G$ is the set of points $\{x \in \R^{n \times n} \ : \ Ax = \textbf{1}, x \geq 0\}$ (where \textbf{1} is the all 1s vector).  By the main theorem of Trubin \cite{trubin} and independently \cite{balaspadberg72}, polytopes given by such systems are quasi-integral  (see also Theorem 7.2 in Chapter 4 of \cite{Yemelichev}).  Therefore, we need to rewrite the defining equations presented in Proposition~\ref{AutGEncodingprop}
to fit this desired shape.
Fix indices $1 \leq i,j \leq n$ and consider the row of $P_G$ defined by the equation $$\sum_{r \in \delta(j)} 
P_{ir} - \sum_{k \in \delta(i)} P_{kj} = 0.$$  
Here $\delta(i)$ denotes those vertices $j$ which are connected to $i$.
Adding the equation $\sum_{r=1}^{n} P_{rj} = 1$ to both sides of this expression yields 
\begin{equation}\label{newProp41eqs}
\sum_{r \in \delta(j)} P_{ir} +\sum_{k \notin \delta(i)} P_{kj} = 1.  
\end{equation}

We can therefore replace the original $n^2$ equations defining $P_G$ by (\ref{newProp41eqs}) over all $1 \leq i,j \leq n$.  The result now follows provided that no summand in each of these equations repeats.  
However, this is clear since if summands $P_{kj}$ and $P_{ir}$ are the same, then $r=j$, which is impossible since $r \in \delta(j)$.
\end{proof}

We would still like to find a tighter description of $IP_G$ in terms of inequalities. For this purpose, recall the radical polynomial ideal $I_G$ in Proposition \ref{AutGEncodingprop} and its real variety $V_{\R}(I_G)$.
We approximate a tighter description of $IP_G$ using a hierarchy of projected semidefinite  relaxations of $conv(V_{\R}(I_G))$. When these relaxations are tight, we obtain a full description of $IP_G$ that allows us to optimize and determine feasibility via semidefinite programming.  

We begin with some preliminary definitions from \cite{GouveiaParriloThomas2008} and motivated by Lov\'{a}sz \& Schrijver \cite{LovaszSchrijver1991}.
Let $I \subset \mathbb{R}[x_1,\ldots,x_n]$ be a \textit{real radical ideal} ($I = \mathcal{I}(V_{\R}(I))$). 
A polynomial $f$ is said to be \textit{nonnegative} mod $I$ (written $f \geq 0$ (mod $I$)) if $f(p) \geq 0$ for all $p \in V_{\mathbb{R}}(I)$. 
Similarly, a polynomial $f$ is said to be a \textit{sum of squares} mod $I$ if there exist $h_1,\ldots,h_m \in \mathbb{R}[x_1,\ldots,x_n]$ such that $f - \sum_{i=1}^{m} h_{i}^{2} \in I$. 
 If the degrees of the $h_1,\ldots,h_m$ are bounded by some positive integer $k$, we say $f$ is $k$-sos mod $I$.  

 The \emph{$k$-th theta body} of $I$, denoted $TH_k(I)$, is the subset
 of $\R^n$ that is nonnegative on each $f \in I$ that is $k$-sos mod
 $I$.  We say that a real variety $V_{\mathbb{R}}(I)$ is \textit{theta
   $k$-exact} if $\overline{conv(V_{\mathbb{R}}(I))} = TH_k(I)$.
   When the ideal $I$ is real radical, theta bodies provide a hierarchy of
 semidefinite relaxations of $\overline{ conv(V_{\R}(I))}$:

$$TH_1(I) \supseteq TH_2(I) \supseteq \cdots \supseteq \overline{ conv(V_{\R}(I))}$$
because in this case theta bodies can be expressed as projections of feasible regions of
 semidefinite programs (such regions are called \emph{spectrahedra}).  In order to 
 exploit this theory, we must establish that $I_G$ is indeed real radical.
 
 \begin{lemma}\label{realradical}
 The ideal $I_G \subseteq \R[x_1,\ldots,x_n]$ is real radical.
 \end{lemma} 
 \begin{proof}
Let $J_G$ be the ideal in $\C[x_1,\ldots,x_n]$ generated by the same polynomials that generate $I_G$, and $\sqrt[\R]{I_G}$ be the real radical of $I_G$.  Since the polynomial $x_i^2 - x_i \in J_G$ for each $1 \leq i \leq n$, Lemma 2.1 of \cite{HillarWindfeldt08} implies $J_G = \sqrt{J_G}$ (where $\sqrt{J_G}$ is the radical of $J_G$).  Together with the fact that $V_{\C}(J_G) = V_{\R}(I_G)$, this implies $J_G \supseteq \sqrt[\R]{I_G}$.  Since $I_G = J_G \cap \R[x_1,\ldots,x_n]$, we conclude $I_G \supseteq \sqrt[\R]{I_G}$.  The result follows since trivially, $I_G \subseteq \sqrt[\R]{I_G}$.
 \end{proof}
 
From Lemma~\ref{realradical}, we conclude that if $I_G$ is theta $k$-exact, linear optimization over
the automorphisms can be performed using semidefinite programming
provided that one first computes a basis for the quotient ring
$\R[P_{11} , P_{12} ,\dots , P_{nn} ]/I_G$ (e.g., obtained from the
standard monomials of a Gr\"obner basis). Using such a basis one can
set up the necessary semidefinite programs (see Section 2 of
\cite{GouveiaParriloThomas2008} for details).  In fact, for $k$-exact
ideals, one only needs those elements of the basis up to degree $2k$.
This motivates the need for characterizing those graphs for which $I_G$ is
$k$-exact.

In this section we focus on finding graphs $G$ such that $I_G$ is 1-exact; we shall call such graphs \emph{exact} in what follows.   The key to finding exact graphs is the following combinatorial-geometric characterization.  

\begin{theorem}\cite{GouveiaParriloThomas2008} \label{compressed}
Let $V_{\mathbb{R}}(I) \subset \mathbb{R}^n$ be a finite real variety.  Then $V_{\mathbb{R}}(I)$ is exact if and only if there is a finite linear inequality description of $conv(V_{\R}(I))$ such that for every inequality $g(x) \geq 0$, there is a hyperplane $g(x) = \alpha$ such that every point in $V_{\mathbb{R}}(I)$ lies either on the hyperplane $g(x)=0$ or the hyperplane $g(x) = \alpha$.
\end{theorem}

A result of Sullivant (see Theorem 2.4 in \cite{Sullivant2006}) directly implies that when the polytope $P=conv(V_{\R}(I))$ 
is lattice isomorphic to an integral polytope of the form  $[0,1]^n \cap  L$ where $L$ is an affine subspace, then $P$ satisfies 
the condition of Theorem \ref{compressed}. Putting these ideas together we can prove compactness implies exactness.
Furthermore, the class of exact graphs properly extends the class of
compact graphs.  The proof of this latter fact is an extension of a result in \cite{tinhofer}.

\begin{theorem} \label{compactinsideexact}
The class of exact graphs strictly contains the class of compact graphs. More precisely:
\begin{enumerate} 

\item If $G$ is a compact graph, then $G$ is also exact. 

\item Let $G_1,\ldots,G_m$ be $k$-regular connected compact graphs, and let $G = \bigsqcup_{i=1}^m G_i$ be the graph that is the disjoint union of $G_1,\ldots,G_m$.  Then 
$G$ is always exact, but $G$ may not be compact.  Indeed, $G$ is compact if and only if $G_i \cong G_j$ for all $1 \leq i,j \leq n$. 
\end{enumerate}
\end{theorem}

\begin{proof}
If $G$ is compact, then the integer hull of $P_G$ is precisely the affine space $$\{P \in \R^{n \times n} \ : \ P A_G = A_G P, \ \sum_{i=1}^{n} P_{ij} = \sum_{j=1}^{n} P_{ij} = 1, \ 1 \leq i,j \leq n \}$$ intersected with the cube ${[0,1]}^{n \times n}$.  That $G$ is exact  follows from Theorem 2.4 of \cite{Sullivant2006}. 

We now prove Statement  2. If $G_i \not\cong G_j$ for some pair $(i,j)$, then $G$ was shown to be non-compact by Tinhofer (see \cite[Lemma 2]{tinhofer}).  Nevertheless, $G$ is exact. We prove this for $m=2$, and the result will follow by induction.  We claim that if $G = G_1 \sqcup G_2$ with $G_1 \not\cong G_2$, then the integer hull $IP_G$ is the solution set to the following system (which we denote by $\tilde{IP}_G$): 
\begin{eqnarray*}~\label{inthull}
{(PA_G - A_GP)}_{i,j} = 0& \ \ \ &1 \leq i,j \leq n, \notag \\
\sum_{i=1}^n P_{i,j} = 1& \ \ \ &1 \leq j \leq n, \\
\sum_{j=1}^n P_{i,j} = 1& \ \ \ &1 \leq i \leq n, \\
\sum_{i=1}^{n_1} \sum_{j=n_1+1}^{n_1+n_2} P_{i,j} = 0, \\
0 \leq P_{i,j} \leq 1,
\end{eqnarray*}
where $n_i=|V(G_i)|$ with $n_1 \leq n_2$. Statement 2 then follows again from Theorem 2.4 of \cite{Sullivant2006}.    

We now prove the claim.  Let $A_{G_i}$ be the adjacency matrix of $G_i$.  
Index the adjacency matrix of $G = G_1 \sqcup G_2$ so that the first $n_1$ rows (and hence first
$n_1$ columns) index the vertices of $G_1$.
Any feasible $P$ of $P_G$ can be written as a block matrix 
\[P = \begin{pmatrix} A_P & B_P \\ C_P & D_P \end{pmatrix},\] 
in which $A_P$ is $n_1 \times n_1$.  Since $G_1$ and $G_2$ are not isomorphic, the only
integer vertices of $P_G$ are of the form $\begin{pmatrix} P_1 & 0 \\
0 & P_2 \end{pmatrix}$ where $P_i$ is an automorphism of $G_i$.  

Now let $P$ be any non-integer vertex of $P_G$.  We claim that the row sums of $B_P$ must be 1.  This will establish that $IP_G$ is described by the system $\tilde{IP}_G$.  To see this, observe that if $Q$ is any point in $P_G$ not in $IP_G$, it is a convex combination of points in $P_G$, one of which (say $P$) is non-integer.  If the row sums of $B_P$ are 1, then $Q$ violates the system $\tilde{IP}_G$.

We now prove that if $P$ is a non-integer vertex of $P_G$, then the row sums of $B_P$ must be 1.  Since $P$ commutes with the adjacency matrix $A_G$ of $G$, we must have
$$ A_P A_{G_1} = A_{G_1} A_P, \ \ B_P A_{G_2} = A_{G_1} B_P, \ \ 
C_P A_{G_2} = A_{G_1} C_P, \ \ D_P A_{G_2} = A_{G_2} D_P. $$
Let $\{b_1,\ldots,b_{n_2} \}$ be the column sums of $B_P$.  We shall calculate the sum of  the entries in each column of $B_PA_{G_2} = A_{G_1}B_P$ in two ways.  First, consider $A_{G_1}B_P$.  Since $G_1$ is $k$-regular, each entry of the $i$-th column of $B_P$ will contribute exactly $k$ times to the sum of the entries of the  $i$-th column of 
$A_{G_1}B_P$.  Thus, the sum of the entries of the $i$-th column of $A_{G_1}B_P$ is $kb_i$.  

Second, consider $B_PA_{G_2}$.  The sum of the entries in its $i$-th column is the sum of 
the entries of the columns of $B_P$ indexed by the neighbors of $i$ in $G_2$.  Thus, the sum of the entries in the $i$-th column of $B_PA_{G_2}$ is $\sum_{l \in 
\delta_{G_2}(i)} b_l$.  It follows that $k b_i = \sum_{l \in \delta_{G_2}(i)} b_l$ for each $1 \leq i \leq n$.  This equality can be written concisely as:
$$\begin{pmatrix}
kI_{n_2 \times n_2} - A_{G_2}
\end{pmatrix}
\begin{pmatrix}
b_1 \\ \vdots \\ b_{n_2}
\end{pmatrix}
= 0.
$$ 
The matrix $kI_{n_2 \times n_2} - A_{G_2}$ is the Laplacian of $G_2$. It is well known that the kernel of the
Laplacian of a connected graph is one dimensional (see \cite{godsil}, Lemma 13.1.1).  Since $G_2$ is regular, the kernel contains the all ones vector.  It follows that $b_1 = \cdots = b_{n_2}$.  By a similar argument, the row sums of $C_P$ are all the
same.  Since all row sums and column sums of $P$ are 1, and the row
sums and column sums of $A_{G_1}$ are the same, it follows that the row
sums of $B_P$ are equal and are the same as the column sums of $C_P$.  

Now assume for contradiction that the row sums of $B_P$ are not 1.  If the row sums are 0,
then $B_P$ and $C_P$ would be 0 matrices.  Since $G_1$ and $G_2$ are compact
this would imply $A_P$ and $D_P$ are permutation matrices, contradicting that $P$ is not integral.
Thus the sum of each row of $B_P$ is $\lambda$ with $0 < \lambda <
1.$ This implies the sum of the rows of $A_P$ is $1 - \lambda$ and that $\frac{1}{1 - \lambda} A_P$ is a feasible solution to $P_{G_1}$.  By
compactness of $G_1$, the matrix $\frac{1}{1 - \lambda} A_P$ is a convex
combination $\sum_{i=1}^{k} \mu_k Q_k$ of permutations $Q_k$ of $G_1$. This implies that
$$P = \sum_{i=1}^{k} \mu_i \begin{pmatrix} (1 - \lambda)Q_k & B_P \\ C_P & D_P  \end{pmatrix},$$ 
which is a convex combination of feasible solutions to $P_G$, contradicting $P$ being a vertex.  It follows that the row sums of $B_P$ must be 1.
\end{proof}

Exact graphs are then more abundant than compact graphs and the convex hull of automorphisms of an exact graph has
a description in terms of semidefinite programming. It is thus desirable to find nice classes of graphs that are exact. 
Notice that exactness is really a property of the set of permutation matrices representing an automorphism group. 
This discussion motivates the following question.  
\begin{question}
Which permutation subgroups of $S_n$ are exact?
\end{question}

Here we view a permutation subgroup of $S_n$ through its natural permutation representation in $\R^{n \times n}$.  
In this light, a permutation subgroup can be considered as a variety, and we say the permutation subgroup is \textit{exact}
if this variety is exact. As an example, consider the alternating group $A_n$ as a subgroup of $S_n$.  It is known (see \cite{cameron}) that 
$A_n$ is never the automorphism group of a graph on $n$ vertices, so it cannot be presented as the integer points of 
a polytope of the form $P_G$ with $|V(G)|=n$.  However, there is a description of $A_n$ as a variety whose points
are vertices of the $n \times n$ Birkhoff polytope:

\begin{equation*}
\begin{split}
\sum_{j=1}^{n} P_{i,j} &= 1, \  \ 1 \leq i \leq n; \ \ \  \sum_{i=1}^{n} P_{i,j} = 1, \ \ 1 \leq j \leq n; \\
\det(P) &= 1; \ \ \ \ \ \ \ \ \ \ \ \ \ \ \ {P}^{2}_{i,j} - P_{i,j} = 0, \ \ 1 \leq i,j \leq n. \\
\end{split}
\end{equation*}

More generally, when a finite permutation group has a description as a variety, we can apply the theory of theta bodies to obtain descriptions of convex hulls.  Using the algebraic-geometric ideas outlined in \cite{sturmfels} we give a sufficient condition 
for exactness of permutation groups.

Let $A = \{\sigma_1,\ldots,\sigma_d\}$ be a subgroup of $S_n$.  We consider $A$ as the set of matrices $\{P_{\sigma_1},\ldots,P_{\sigma_d}\}$ $\subseteq \mathbb{Z}^{n \times n}$, where $P_{\sigma_i}$ is the permutation matrix corresponding to $\sigma_i$.  Let $\C[\textbf{x}] := \C[x_{\sigma_1},\ldots,x_{\sigma_d}]$ be the polynomial ring in $d$ indeterminates indexed by permutations in $A$, and let $\C[\textbf{t}] := \C[t_{ij}: 1 \leq i,j \leq n]$. 

 The algebra homomorphism induced by the map 
\begin{equation}\label{toricalghom}
\pi : \C[\textbf{x}] \to \C[\textbf{t}] , \  \pi(x_{\sigma_i}) = \prod_{ 1 \leq j,k \leq n}  t_{jk}^{(P_{\sigma_i})_{jk} }
\end{equation}
has kernel $I_A$, which is a prime \textit{toric ideal} \cite{sturmfels}.
By Theorem \ref{compressed}, Corollary 8.9 in \cite{sturmfels}, and Corollary 2.5 in \cite{Sullivant2006},
the group $A$ is exact if and only if for every reverse lexicographic term ordering $\prec$ on $\C[x]$, the initial ideal $in_{\prec}(I_A)$ is generated by square-free monomials.  We now describe a family of permutation groups that are exact.

Let $A \subseteq \mathbb{Z}^{n \times n}$ be a subgroup of $S_n$.  We say that $A$ is \textit{permutation summable} if for any permutations $P_1,\ldots,P_m \in A$ satisfying the inequality $\sum_{i=1}^m P_i - I \geq 0$ (entry-wise), we have that $\sum_{i=1}^{m} P_i -I$ is 
also a sum of permutation matrices in $A$. For example, Birkhoff's Theorem (see  e.g., Theorem 1.1 in Chapter 5 of \cite{Yemelichev})
 implies $S_n$ is  permutation summable. Note that in this case $P_{S_n}$ is the Birkhoff polytope which is known to be exact by the results in \cite{GouveiaParriloThomas2008}.  We prove the following result.

\begin{theorem}\label{groebner} 
Let $A = \{\sigma_1,\ldots,\sigma_d\}$ be a permutation group that is a subgroup of  $S_n$.

(1) If $A$ is permutation summable, then $A$ is exact.

(2) Suppose $I_A$, the toric ideal associated to $A$, has a quadratically generated 
Gr\"obner basis with respect to any reverse lexicographic ordering $\prec$, then $A$ is exact.
\end{theorem}

\begin{proof}
Let $I_A$ be the kernel of the algebra homomorphism induced by (\ref{toricalghom}).  We shall abbreviate the action of $\pi$ on $x_{\sigma}$ by $\pi(x_{\sigma}) = t^{P_\sigma}$ for any $\sigma \in A$. 

 Let $\mathfrak{G}$ be a reduced Gr\"obner basis for $I_A$ with respect to some reverse lexicographic order $\prec$ on $\{x_{\sigma_1},\ldots, x_{\sigma_d}\}$.  Let $x^u - x^v \in \mathfrak{G}$ with leading term $x^u$.  By Theorem \ref{compressed}, Corollary 8.9 in \cite{sturmfels} and Corollary 2.5 in \cite{Sullivant2006}, Statement (1) follows if we can find a square-free monomial $x^{u'} \in in_{\prec}(I_A)$ such that $x^{u'}$ divides $x^u$.  

Let $x_{\tau}$ be the smallest variable dividing $x^v$  with respect to $\prec$.  Then $x_{\tau}$ is smaller than any variable appearing in $x^u$ by the choice of a reverse lexicographic ordering.  Since $x^u - x^v \in \mathfrak{G}$, we have $\pi(x^u) = \pi(x^v)$.  It follows that $\pi(x_{\tau})$ divides $\pi(x^u)$, so letting $x^u = x_{\sigma_{i_1}} \cdots x_{\sigma_{i_k}}$ for some $\{\sigma_{i_1},\ldots,\sigma_{i_k}\} \subseteq A$, we have  
$$\dfrac{\pi(x^u)}{\pi(x_{\tau})} = t^{P_{\sigma_{i_1}} + \cdots + P_{\sigma_{i_k}}}  t^{- P_{\tau}},$$ 
in which $\sum_{j=1}^{k} P_{\sigma_{i_j}} - P_{\tau}$ is a matrix with nonnegative integer entries.  Choose  a subset  $\{\rho_1,\ldots,\rho_r\}$ $\subset \{\sigma_{i_1},\ldots,\sigma_{i_k}\}$ such that $\{P_{\rho_1},\ldots,P_{\rho_r}\}$ minimally supports $P_{\tau}$ with $P_{\rho_i} \neq P_{\rho_j}$ for all $i,j$, and let $x^{u'} = x_{\rho_1} \cdots x_{\rho_r}$.  We claim that $x^{u'}$ is a square-free monomial that divides $x^u$ and lies in $in_{\prec}(I_A)$, which will prove Statement (1).  

By construction, all indeterminates $x_{\rho_1},\ldots,x_{\rho_r}$ are distinct, so $x^{u'}$ is square-free.  Moreover, since $\{\rho_1,\ldots,\rho_r\} \subset \{\sigma_{i_1},\ldots,\sigma_{i_k}\}$, we have that $x^{u'}$ divides $x^u$.    It remains to show that $x^{u'}$ lies in $in_{\prec}(I_A)$.  To see this, note that $\sum_{i=1}^{r} P_{\rho_i} - P_{\tau}$ 
has nonnegative integer entries, and hence so does 
\[M = \sum_{i=1}^{r} {(P_{\tau})}^{-1}P_{\rho_i} - I\] (multiplying by $P_{\tau}^{-1}$ permutes matrix entries, and therefore does not effect nonnegativity). Since $A$ is  permutation summable, the matrix $M$ is a sum of matrices in $A$, and hence so is $P_{\tau}M$ = $\sum_{i=1}^{r} P_{\rho_i} - P_{\tau}$.  It follows that 
$$\sum_{i=1}^{r} P_{\rho_i} - P_{\tau} = \sum_{j=1}^{r-1}  P_{\sigma_{l_j}}$$ 
for some $\{\sigma_{l_1},\ldots,\sigma_{l_{r-1}} \} \subset A$.  In particular, $\pi(x^{u'}) = \pi(x_{\tau}) \cdot \pi(x^{v'})$ and so $x^{u'} - x_{\tau} x^{v'} \in I_A$.  Since $x_{\tau}$ is smaller than any term in $x^{u'}$ (the monomial $x^{u'}$ divides $x^u$ and the same holds for $x^u$), the leading term of $x^{u'} - x_{\tau} x^{v'}$ is $x^{u'}$; hence, $x^{u'} \in in_{\prec}(I_A)$.  This proves Statement (1).

For Statement (2), since any Gr\"obner basis is quadratically generated, by part (1) it suffices to show that if $P_1,P_2,Q \in A$ with all entries of $P_1 + P_2 - Q$ nonnegative, then $P_1 + P_2 - Q$ is a permutation matrix.  Since $supp(Q) \subset supp(P_1) \cup supp(P_2)$, the permutation $Q$ is a vertex of a face containing $P_1$ and $P_2$.  By Theorem 3.5 of \cite{guralnick}, $Q$ is on the smallest face containing $P_1$ and $P_2$, and this face is centrally symmetric.  Thus, there is a vertex $R$ such that $Q + R = P_1 + P_2$, and the result follows.

\end{proof}

In light of Theorem~\ref{groebner}, we would like to find permutation groups $A$ that are permutation summable.  As we have seen, Birkhoff's Theorem (see \cite{sturmfels}) implies that $S_n$ is permutation summable.  We can use this fact to construct more permutation summable groups.  For instance, $S_{n_1} \times \cdots \times S_{n_m}$ is permutation summable, simply by applying the permutation summability condition on each $S_{n_i}$ and taking direct sums.  More generally, if $H_{1}, \ldots, H_{m}$ are permutation summable, then so is $H_{1} \times \cdots \times H_{m}$.  We present another class of permutation summable groups that contains familiar groups.

\begin{definition} 
Let $A$ be a permutation subgroup of $S_n$.  We say $A$ is \textit{strongly fixed-point free} if for every $\sigma \in A \backslash \{1\}$, we have $\sigma(i) \neq i$ for any $i \in \{1,\ldots,n\}$.
\end{definition}

\begin{corollary}
Let $A$ be a strongly fixed-point free subgroup of $S_n$.  Then $A$ is exact.
\end{corollary}

\begin{proof} Let $A$ be strongly fixed-point free. 
Consider any subset $\{P_{\sigma_1},\ldots,P_{\sigma_k}\}$ of $A$ and assume $\sum_{i=1}^{k} P_{\sigma_i} - I$ is a matrix with nonnegative entries.  Then one of the matrices in $A$ contains a fixed point.  Without loss of generality, assume $P_{\sigma_1}$ is one such matrix.  Since $A$ is strongly fixed-point free, we have $P_{\sigma_1} = I$.  Hence,
$$\sum_{i=1}^{k} P_{\sigma_i} - I = \sum_{i=2}^{k} P_{\sigma_i},$$ and thus $A$ is permutation summable.  The result now follows from Theorem~\ref{groebner}.
\end{proof}

There are many well-known families of permutation groups that are strongly fixed-point free, and hence exact.  These include the group generated by any $n$ cycle in $S_n$, and even dihedral groups (dihedral groups of order $4n$ as subgroups of $S_{2n}$).

\section{Acknowledgements}
We would like to thank the referee for his or her valuable suggestions and corrections.

\bibliography{references}{}

\providecommand{\bysame}{\leavevmode\hbox to3em{\hrulefill}\thinspace}
\providecommand{\MR}{\relax\ifhmode\unskip\space\fi MR }
\providecommand{\MRhref}[2]{%
  \href{http://www.ams.org/mathscinet-getitem?mr=#1}{#2}
}
\providecommand{\href}[2]{#2}
\begin{thebibliography}{10}

\bibitem{alonsurvey}
N.~Alon, \emph{Combinatorial {N}ullstellensatz}, Combin. Probab. Comput.
  \textbf{8} (1999), no.~1-2, 7--29, Recent trends in combinatorics
  (M\'atrah\'aza, 1995).

\bibitem{AlonTarsi92}
N.~Alon and M.~Tarsi, \emph{Colorings and orientations of graphs},
  Combinatorica \textbf{12} (1992), no.~2, 125--134.

\bibitem{avis}
D.~Avis and K.~Fukuda, \emph{Reverse search for enumeration}, Discrete Appl.
  Math. \textbf{65} (1996), no.~1-3, 21--46, First International Colloquium on
  Graphs and Optimization (GOI), 1992 (Grimentz).

\bibitem{balaspadberg72}
E.~Balas and M.W. Padberg, \emph{On the set-covering problem}, Operations Res.
  \textbf{20} (1972), 1152--1161.

\bibitem{Bayer82}
D.A. Bayer, \emph{The {D}ivision {A}lgorithm and the {H}ilbert {S}cheme}, Ph.D.
  thesis, Harvard University, 1982.

\bibitem{Cameron2001}
K.~Cameron, \emph{Thomason's algorithm for finding a second hamiltonian circuit
  through a given edge in a cubic graph is exponential on krawczyk's graphs},
  Discrete Math. \textbf{235} (2001), no.~1-3, 69--77, Combinatorics (Prague,
  1998).

\bibitem{cameron}
P.J. Cameron, \emph{Automorphisms of graphs}, Topics in Algebraic Graph Theory
  (R.J.~Wilson L.W.~Beineke, ed.), Cambridge Univ. Press, 2004, pp.~203--221.

\bibitem{godsil}
A.~Chan and C.~Godsil, \emph{Graph symmetry: Algebraic methods and
  applications}, ch.~4, pp.~75--106, Kluwer Academic Publishers, Montr\'{e}al,
  QC, Canada., 1997.

\bibitem{CleggEdmondsImpagliazzo1996}
M.~Clegg, J.~Edmonds, and R.~Impagliazzo, \emph{Using the {G}roebner basis
  algorithm to find proofs of unsatisfiability}, STOC '96: Proceedings of the
  twenty-eighth annual ACM symposium on Theory of computing (New York, NY,
  USA), ACM, 1996, pp.~174--183.

\bibitem{CourtoisKlimovPatarinShamir2000}
N.~Courtois, A.~Klimov, J.~Patarin, and A.~Shamir, \emph{Efficient algorithms
  for solving overdefined systems of multivariate polynomial equations},
  Advances in cryptology---EUROCRYPT 2000 (Bruges) (Berlin), Lecture Notes in
  Comput. Sci., vol. 1807, Springer, 2000, pp.~392--407.

\bibitem{CoxLittleOShea92}
D.~Cox, J.~Little, and D.~O'Shea, \emph{Ideals, varieties, and algorithms}, 3
  ed., Undergraduate Texts in Mathematics, Springer, 2007, An introduction to
  computational algebraic geometry and commutative algebra.

\bibitem{DeLoera95}
J.~A. De~Loera, \emph{Gr\"obner bases and graph colorings}, Beitr\"age Algebra
  Geom. \textbf{36} (1995), no.~1, 89--96.

\bibitem{DeLoeraLeeMalkinMargulies08}
J.A. De~Loera, J.~Lee, P.N. Malkin, and S.~Margulies, \emph{Hilbert's
  {N}ullstellensatz and an algorithm for proving combinatorial infeasibility},
  Proceedings of the Twenty-first International Symposium on Symbolic and
  Algebraic Computation (ISSAC 2008), 2008.

\bibitem{DeLoeraLeeMarguliesOnn08}
J.A. De~Loera, J.~Lee, S.~Margulies, and S.~Onn, \emph{Expressing combinatorial
  problems by systems of polynomial equations and {H}ilbert's
  {N}ullstellensatz}, Combin. Probab. Comput. \textbf{18} (2009), no.~4,
  551--582.

\bibitem{ima}
J.A. De~Loera, P.~Malkin, and P.~Parrilo, \emph{Computation with polynomial
  equations and inequalities arising in combinatorial optimization},
  http://arxiv.org/abs/0909.0808, 2009.

\bibitem{dummit}
D.S. Dummit and R.~M. Foote, \emph{Abstract algebra}, 3 ed., John Wiley \& Sons
  Inc., 2004.

\bibitem{Eliahou92}
S.~Eliahou, \emph{An algebraic criterion for a graph to be four-colourable},
  International {S}eminar on {A}lgebra and its {A}pplications ({S}panish)
  ({M}\'exico {C}ity, 1991), Aportaciones Mat. Notas Investigaci\'on, vol.~6,
  Soc. Mat. Mexicana, M\'exico, 1992, pp.~3--27.

\bibitem{Faugere1999}
J.~C. Faug{\'e}re, \emph{A new efficient algorithm for computing {G}r\"obner
  bases $({F}\sb 4)$}, J. Pure Appl. Algebra \textbf{139} (1999), no.~1-3,
  61--88, Effective methods in algebraic geometry (Saint-Malo, 1998).

\bibitem{Fischer88}
K.G. Fischer, \emph{Symmetric polynomials and {H}all's theorem}, Discrete Math.
  \textbf{69} (1988), no.~3, 225--234.

\bibitem{friedlander}
S.~Friedland, \emph{Graph isomorphism and volumes of convex bodies},
  http://arxiv.org:0911.1739, 2009.

\bibitem{GouveiaParriloThomas2008}
J.~Gouveia, P.A. Parrilo, and R.R. Thomas, \emph{Theta bodies for polynomial
  ideals}, http://arxiv.org:0809.3480, 2008.

\bibitem{guralnick}
R.~M Guralnick and D.~Perkinson, \emph{Permutation polytopes and indecomposable
  elements in permutation groups}, J. Combin. Theory Ser. A \textbf{113}
  (2006), no.~7, 1243--1256.

\bibitem{HillarLim}
C.~J. Hillar and L-H. Lim, \emph{Most tensor problems are \text{NP} hard},
  preprint, 2010.

\bibitem{HillarWindfeldt08}
C.~J. Hillar and T.~Windfeldt, \emph{Algebraic characterization of uniquely
  vertex colorable graphs}, J. Combin. Theory Ser. B \textbf{98} (2008), no.~2,
  400--414.

\bibitem{KehreinKreuzer05}
A.~Kehrein and M.~Kreuzer, \emph{Characterizations of border bases}, J. Pure
  Appl. Algebra \textbf{196} (2005), no.~2-3, 251--270.

\bibitem{Kollar1988}
J.~Koll{\'a}r, \emph{Sharp effective {N}ullstellensatz}, J. Amer. Math. Soc.
  \textbf{1} (1988), no.~4, 963--975.

\bibitem{Yemelichev}
M.~M. Koval{\"e}v, M.~K. Kravtsov, and V.A. Yemelichev, \emph{Polytopes, graphs
  and optimisation}, Cambridge University Press, Cambridge, 1984, Translated
  from the Russian by G. H. Lawden.

\bibitem{Lasserre2002}
J.~B. Lasserre, \emph{An explicit equivalent positive semidefinite program for
  nonlinear $0$-$1$ programs}, SIAM J. Optim. \textbf{12} (2002), no.~3,
  756--769 (electronic).

\bibitem{Laurent2007}
M.~Laurent, \emph{Semidefinite representations for finite varieties}, Math.
  Program. \textbf{109} (2007), no.~1, Ser. A, 1--26.

\bibitem{moniquefranz}
M.~Laurent and F.~Rendl, \emph{Semidefinite programming \& integer
  programming}, Handbook on Discrete Optimization (K.~Aardal, G.~Nemhauser, and
  R.~Weismantel, eds.), Elsevier B.V., 2005, pp.~393--514.

\bibitem{LiLi81}
S.-Y.R. Li and W.C.W Li, \emph{Independence numbers of graphs and generators of
  ideals}, Combinatorica \textbf{1} (1981), no.~1, 55--61.

\bibitem{Lovasz1994}
L.~Lov{\'a}sz, \emph{Stable sets and polynomials}, Discrete Math. \textbf{124}
  (1994), no.~1-3, 137--153, Graphs and combinatorics (Qawra, 1990).

\bibitem{LovaszSchrijver1991}
L.~Lov{\'a}sz and A.~Schrijver, \emph{Cones of matrices and set-functions and
  $0$-$1$ optimization}, SIAM J. Optim. \textbf{1} (1991), no.~2, 166--190.

\bibitem{SusanPhd}
S.~Margulies, \emph{Computer algebra, combinatorics, and complexity: Hilbert's
  {N}ullstellensatz and {NP}-complete problems}, Ph.D. thesis, UC Davis, 2008.

\bibitem{matiyasevich1}
Y.~Matiyasevich, \emph{A criteria for colorability of vertices stated in terms
  of edge orientations}, Discrete Analysis \textbf{26} (1974), 65--71.

\bibitem{matiyasevich2}
\bysame, \emph{Some algebraic methods for calculation of the number of
  colorings of a graph}, Zapiski Nauchnykh Seminarov POMI \textbf{293} (2001),
  193--205.

\bibitem{MourrainTrebuchet08}
B.~Mourrain and P.~Tr{\'e}buchet, \emph{Stable normal forms for polynomial
  system solving}, Theoret. Comput. Sci. \textbf{409} (2008), no.~2, 229--240.

\bibitem{Onn}
S.~Onn, \emph{Nowhere-zero flow polynomials}, Journal of Combinatorial Theory,
  Series A \textbf{108} (2004), 205--215.

\bibitem{Parrilo2003}
P.~A. Parrilo, \emph{Semidefinite programming relaxations for semialgebraic
  problems}, Math. Program. \textbf{96} (2003), no.~2, Ser. B, 293--320,
  Algebraic and geometric methods in discrete optimization.

\bibitem{Parrilo2002}
P.A. Parrilo, \emph{An explicit construction of distinguished representations
  of polynomials nonnegative over finite sets}, IfA AUT02-02, ETH Z\"urich,
  2002.

\bibitem{PokuttaSchulz2009}
S.~Pokutta and A.S. Schulz, \emph{On the connection of the {S}herali-{A}dams
  {C}losure and {B}order {B}ases}, 2009, Working Paper, Technische
  Universit\"at Darmstadt / Massachusetts Institute of Technology.

\bibitem{Seymour1979}
P.~D. Seymour, \emph{Sums of circuits}, Graph Theory and Related Topics (1979),
  341--355.

\bibitem{Simis94}
A.~Simis, W.V. Vasconcelos, and R.H. Villarreal, \emph{On the ideal theory of
  graphs}, J. Algebra \textbf{167} (1994), no.~2, 389--416.

\bibitem{Stetter04}
H.~J. Stetter, \emph{Numerical polynomial algebra}, Society for Industrial and
  Applied Mathematics (SIAM), 2004.

\bibitem{sturmfels}
B.~Sturmfels, \emph{Gr\"obner bases and convex polytopes}, University Lecture
  Series, vol.~8, American Mathematical Society, 1996.

\bibitem{Sullivant2006}
S.~Sullivant, \emph{Compressed polytopes and statistical disclosure
  limitation}, Tohoku Math. J. (2) \textbf{58} (2006), no.~3, 433--445.

\bibitem{Szekeres1973}
G.~Szekeres, \emph{Polyhedral decompositions of cubic graphs}, Bull. Austral.
  Math. Soc. \textbf{8} (1973), 367--387.

\bibitem{tinhofer}
G.~Tinhofer, \emph{Graph isomorphism and theorems of birkhoff type}, Computing
  \textbf{36} (1986), 285--300.

\bibitem{trubin}
V.~A. Trubin, \emph{A method of solution of a special form of integer linear
  programming problems}, Dokl. Akad. Nauk SSSR \textbf{189} (1969), 952--954.

\bibitem{Tutte1946}
W.~T. Tutte, \emph{On hamiltonian circuits}, J. London Math. Soc. \textbf{21}
  (1946), 98--101.

\bibitem{Yap2000}
Chee~Keng Yap, \emph{Fundamental problems of algorithmic algebra}, Oxford
  University Press, New York, 2000. \MR{MR1740761 (2000m:12014)}

\end{thebibliography}
\bibliographystyle{amsplain}

\end{document}